\documentclass[12pt]{article}
\usepackage{amssymb}\usepackage{bm}\usepackage{graphicx}
\begin{document}
\title{Solutions to Monge-Kantorovich equations\\ as stationary points of a dynamical system}

\author{Leonid Prigozhin\\Department of Solar Energy and Environmental
Physics\\ Blaustein Institute for Desert Research\\ Ben-Gurion
University, Sede-Boqer Campus 84990, Israel.}

\maketitle


\begin{abstract} Solutions to Monge-Kantorovich
equations, expressing optimality condition in mass transportation
problem with cost equal to distance, are stationary points of a
critical-slope model for sand surface evolution. Using a dual
variational formulation of sand model, we compute both the optimal
transport density and Kantorovich potential as $t\rightarrow
\infty$ limit of evolving sand flux and sand surface,
respectively.
\end{abstract}
\section{Introduction}
The Monge-Kantorovich equations,
\begin{eqnarray}
f^+-f^-=-\nabla\cdot(a\nabla u),\label{m_bal}\\
|\nabla u|\le 1,\ \ \ a\ge 0,\ \ \ |\nabla u|< 1\longrightarrow a=
0, \label{const_rel}
\end{eqnarray}
appear as a condition of optimality in the classical
Monge-Kantorovich problem of optimal mass transportation with the
cost function equal to the distance of transportation
\cite{A,EG,FM}. Here $f^{\pm}$ are two given distributions of
mass,  $a$ is the transport density, and $u$ is the Kantorovich
potential. Let ${\cal M}(R^n)$ be the set of bounded Radon
measures on $R^n$, ${\cal M}_+(R^n)$ the subset of nonnegative
measures, and $\mbox{Lip}_1(R^n)$ the set of real valued Lipschitz
functions with Lipschitz constant not greater than 1. If
$f^{\pm}\in { \cal M}_+(R^n)$ and satisfy $f^+({R}^n)=f^-({R}^n)$,
there exist a transport density and a corresponding potential,
$\{a,u\}\in {\cal M}_+(R^n)\times \mbox{Lip}_1(R^n)$, satisfying
(\ref{m_bal})-(\ref{const_rel}) in a weak sense, see \cite{A}. If
either $f^+$ or $f^-$ are absolutely continuous with respect to
Lebesgue measure in $R^n$ the transport density $a$ is unique and
also absolutely continuous \cite{A}. Obviously, the Kantorovich
potential $u$ is determined up to an additive constant. If the
constant is fixed, the potential is still not unique outside
spt$(a)$. Note that the system (\ref{m_bal})-(\ref{const_rel})
arises also in a mass optimization problem \cite{BB}.

A critical-slope model for sand surface evolution (\cite{P_PRE},
see also \cite{P_EJAM,AEW})
\begin{eqnarray}
f-\partial_t u=-\nabla\cdot(a\nabla u),\ \ \ u|_{t=0}=u_0(x),\label{h_eq}\\
|\nabla u|\le k,\ \ \ a\ge 0,\ \ \ |\nabla u|< k\longrightarrow a=
0, \label{h_rel}
\end{eqnarray}
can be regarded as a non-stationary version of the
Monge-Kantorovich equations. Here $x\in R^2$, $u(x,t)$ is the
evolving sand surface, $f(x,t)$ is the intensity of external
source, $k$ is the tangent of the sand angle of repose, $a(x,t)$
is the sand transport density (if $k=1$) or proportional to it,
and the initial surface $u_0$ satisfies the equilibrium constraint
$|\nabla u_0|\le k$.

Indeed, it can be proved that if the source intensity $f$ is
nonnegative, the term $\partial_t u$ is nonnegative too and, as
was noted in \cite{E_survey}, comparison of the two systems shows
that at each time moment sand from the source is instantly
incorporated into the bulk in a way that minimizes sand transport
under the constraint $|\nabla u|\le k$. The height function $u$
plays here the role of Kantorovich potential.

In this note we explore an opposite approach to this similarity.
Let the source intensity be a time-independent Radon measure,
$f=f^+-f^-$, such that $f^{\pm}\in {\cal M}_+(R^n)$ and satisfy
$f^+({R}^n)=f^-({R}^n)$. Starting with an admissible initial
condition $u_0$, we seek a solution to Monge-Kantorovich equations
as the $t\rightarrow\infty$ limit of solution to
(\ref{h_eq})-(\ref{h_rel}). To solve this latter evolutionary
problem we employ the numerical algorithm \cite{BP} based on a
dual variational formulation of the sand model written for flux of
sand pouring down the evolving sand surface \cite{P_PhysD}.
\section{Variational formulations of sand model}
To study the evolution of sand surface, it is convenient to
exclude from equations (\ref{h_eq})-(\ref{h_rel}) the transport
density $a$, the Lagrange multiplier related to the constraint
$|\nabla u|\le k$, and rewrite the model as a variational
inequality for $u$ alone. For bounded domains and appropriate
boundary conditions, the inequality was obtained in
\cite{P_PRE,P_EJAM}. Similar formulation  of the initial value
problem has been independently derived and studied in \cite{AEW}:
\begin{equation}\begin{array}{c}u(.,t)\in K:\ (\partial_t u-f,\varphi-u)\ge
0,\ \ \forall \varphi\in
K,\\u|_{t=0}=u_0,\end{array}\label{vi1}\end{equation} where
$K=\{\varphi\in L^2(R^2)\ :\ |\nabla \varphi|\le k\
\mbox{a.e.}\}$. We assume both $u_0\in K$ and $f$ have compact
supports. In this case spt$(u)$ remains compact for each $t>0$,
see \cite{AEW}.

The inequality (\ref{vi1}) can be used for efficient computation
of the evolving sand surface $u$. However, the dual variable, $a$,
remains unknown and difficult to find. Such situation is typical
also of other critical-state problems. Because of this reason the
dual variational formulations \cite{P_PhysD}, resembling mixed
variational inequalities in elasto-plasticity \cite{HR} and
allowing to find both variables simultaneously (see \cite{BP}),
can be preferable.
Having in mind the application to Monge-Kantorovich problem, we
assume $x\in R^n$ and allow for non-constant coefficients
$k=k(x)\ge 0$ to treat problems with, say, subregions through
which mass transportation is forbidden or, on the contrary, where
it is free of charge.

Suppose the pair $\{u,a\}$ satisfies the model relations
(\ref{h_eq})-(\ref{h_rel}). Then, for the horizontal projection of
surface flux of sand $\bm{q}=-a\nabla u$ and arbitrary test field
$\bm{\psi}$ we obtain
$$\nabla u\cdot(\bm{\psi}-\bm{q})\ge -|\nabla
u||\bm{\psi}|-\nabla u\cdot\bm{q}=-|\nabla
u||\bm{\psi}|+k|{\bm{q}}|\ge -k|\bm{\psi}|+k|{\bm{q}}|.$$ Noting
that support of $u$ is bounded and integrating, we obtain $(\nabla
u,\bm{\psi}-\bm{q})\ge\phi(\bm{q})-\phi(\bm{\psi}),$ where
$\phi(\bm{q})=\int_{R^n}k|{\bm{q}}|$, and also $(\nabla
u,\bm{\psi}-\bm{q})= -(u,\nabla\cdot\{\bm{\psi}-\bm{q}\})$.
Therefore,
$$\phi(\bm{\psi})-\phi(\bm{q})-(u,\nabla\cdot\{\bm{\psi}-\bm{q}\})
\ge 0.$$ Integrating now the balance equation (\ref{h_eq}) in time
we get
\begin{equation}u=u_0+\int_0^tf\,dt-\nabla\cdot \int_0^t\bm{q}\,dt \label{u_eq}
\end{equation} and arrive at the mixed variational inequality
written for the sand flux alone:
\begin{equation}\bm{q}(.,t)\in V:\ \left(\nabla\cdot\! \int_0^t\bm{q}\,dt-\mathcal{F},
\nabla\cdot\{\bm{\psi}-\bm{q}\}\right)
+\phi(\bm{\psi})-\phi(\bm{q})\ge 0 \label{vi2}\end{equation} for
any $\bm{\psi} \in V$. Here $\mathcal{F}=u_0+\int_0^tf\, dt$ and
we define $$V=\{\bm{\psi} \in [{\mathcal M}(R^n)]^n\ :\
\nabla\cdot\bm{\psi}\in L^2(R^n)\}.$$ Several comments about the
variational formulation (\ref{vi2}) can be necessary.

 The problem is
not coercive in the usual Sobolev spaces. Nevertheless, under
suitable conditions on the source function $f$, existence of a
solution can be established \cite{BP2}.

Since the divergence of flux $\bm{q}$ belongs to $L^2$, we expect
the continuity of normal component of sand flux; the tangential
component of $\bm{q}$ can be discontinuous (this corresponds to
continuity of transport density along the transport rays). To
solve such a problem numerically, one should use the
divergence-conforming finite elements that do not enforce any
additional smoothness on the solution.

Suppose the flux has been found as a solution to (\ref{vi2}).
Obviously, the transport density is easily determined for $k>0$ as
$a=|\bm{q}|/k$. The surface $u$ (potential in the
Monge-Kantorovich problem) can  be calculated by means of the
equation (\ref{u_eq}).
\section{Solution of Monge-Kantorovich equations}
To approximate the inequality (\ref{vi2}) numerically, we smoothed
the non-differen\-ti\-able functional $\phi$ by introducing
$|\bm{q}|_{\varepsilon}=(|\bm{q}|^2+\varepsilon^2)^{1/2}$,
discretized the regularized equality problem in time, employed the
divergence-conforming Raviart-Thomas finite elements of lowest
order (see, e.g., \cite{BC}) and used vertex sampling  on the
nonlinear term. Resulting nonlinear algebraic systems were solved
at each time level iteratively using a form of successive
over-relaxation. We refer to \cite{BP2} for  details and make only
several remarks related to efficient implementation of this
algorithm to Monge-Kantorovich equations.

i) The needed solutions are obtained as the $t\rightarrow\infty$
limit of solutions to (\ref{vi2}) with the time-independent
sources satisfying $f(R^n)=0$. To find this limit, we need not
solve the nonlinear equations at each time level with high
accuracy: only a few iterations are needed before the transition
to a new time level. If, however, the Kantorovich potential is
also of interest by some reason, the time step and the number of
iterations should ensure accurate calculation of the integral in
(\ref{u_eq}).

 ii) Suppose that transportation through an open domain
 $\Omega$ is forbidden. To model this
situation, it is possible to solve the problem in
$R^n\setminus\Omega$ with $\bm{q}_n|_{\partial \Omega}=0$. Another
possibility (used in this work) is to solve also in $\Omega$ but
set $k|_{\Omega}=\infty$. Then, if $\bm{q}$ solves (\ref{vi2}),
the flux $\bm{q}|_{\Omega}$ must be zero. Since the normal
component of flux is continuous, only tangential to domain
boundary flux of mass is permitted: the optimal mass
transportation rounds the obstacle.

iii) Numerical solution of (\ref{vi2}) is possible only in a
bounded domain and it is desirable to make the computational
domain as small as possible. However, if the domain is not large
enough, at some (depending on $u_0$) time moment
spt$(\bm{q}(.,t))$ can reach the domain boundary. To avoid any
material loss through the boundary, we can surround the domain by
an obstacle (see previous comment). On the other hand, we do not
want the artificial obstacle to affect our solution. Since the
transport set in Monge-Kantorovich problem (without obstacles)
consists of straight transport rays that begin in spt$(f^+)$ and
end in spt$(f^-)$,  we can choose any computational domain
containing the convex hull of
spt$(f)=$spt$(f^+)\bigcup\,$spt$(f^-)$ and surround it by an
obstacle. Our choice of domain (and $u_0$) can influence only the
Kantorovich potential outside the establishing in
$t\rightarrow\infty$ limit transport set, where this potential is
not unique.

\begin{figure}[hb]
\begin{center}
\includegraphics[width=7.0cm,height=7.0cm]{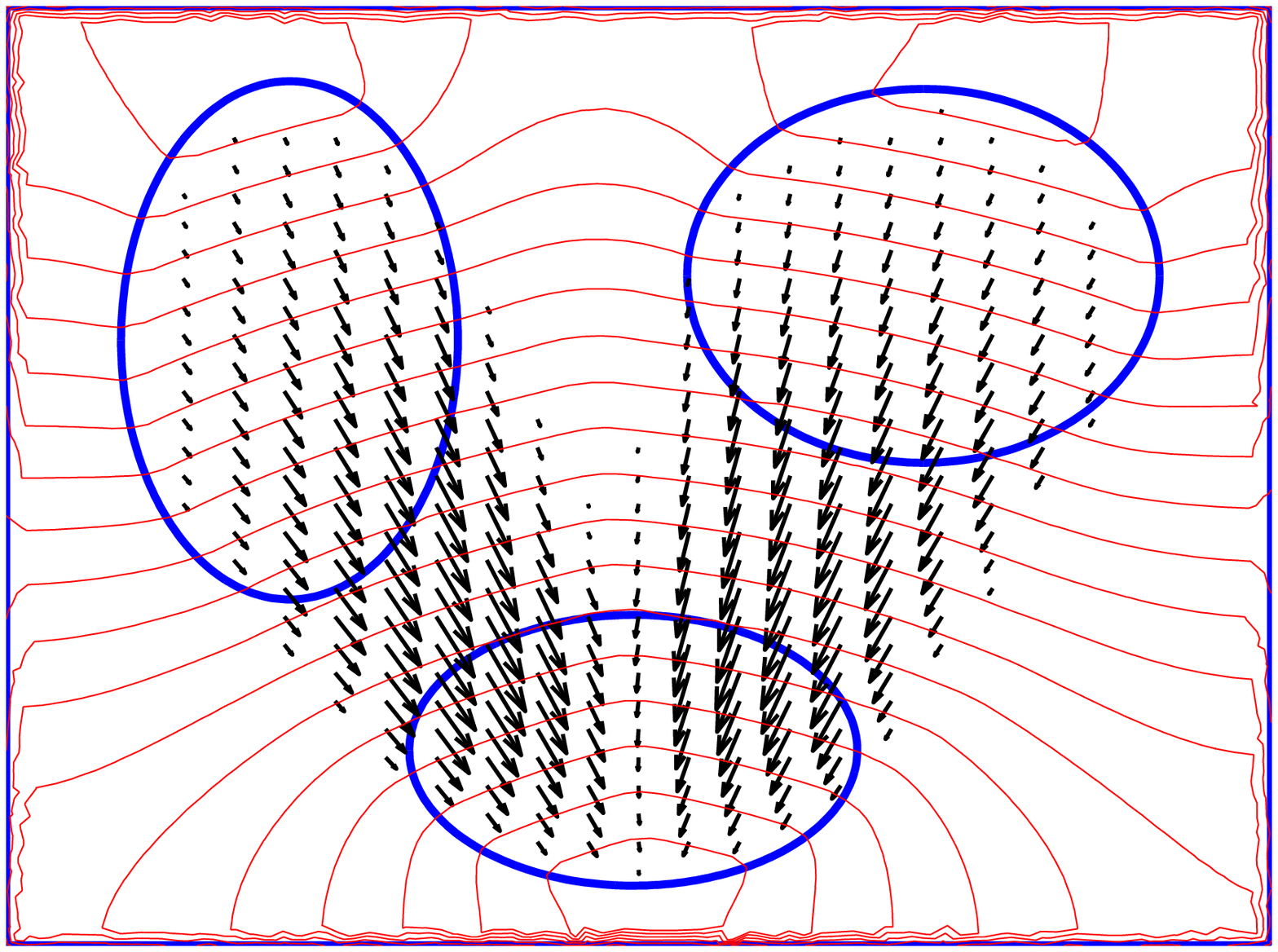}\includegraphics[width=7.0cm,height=7.0cm]{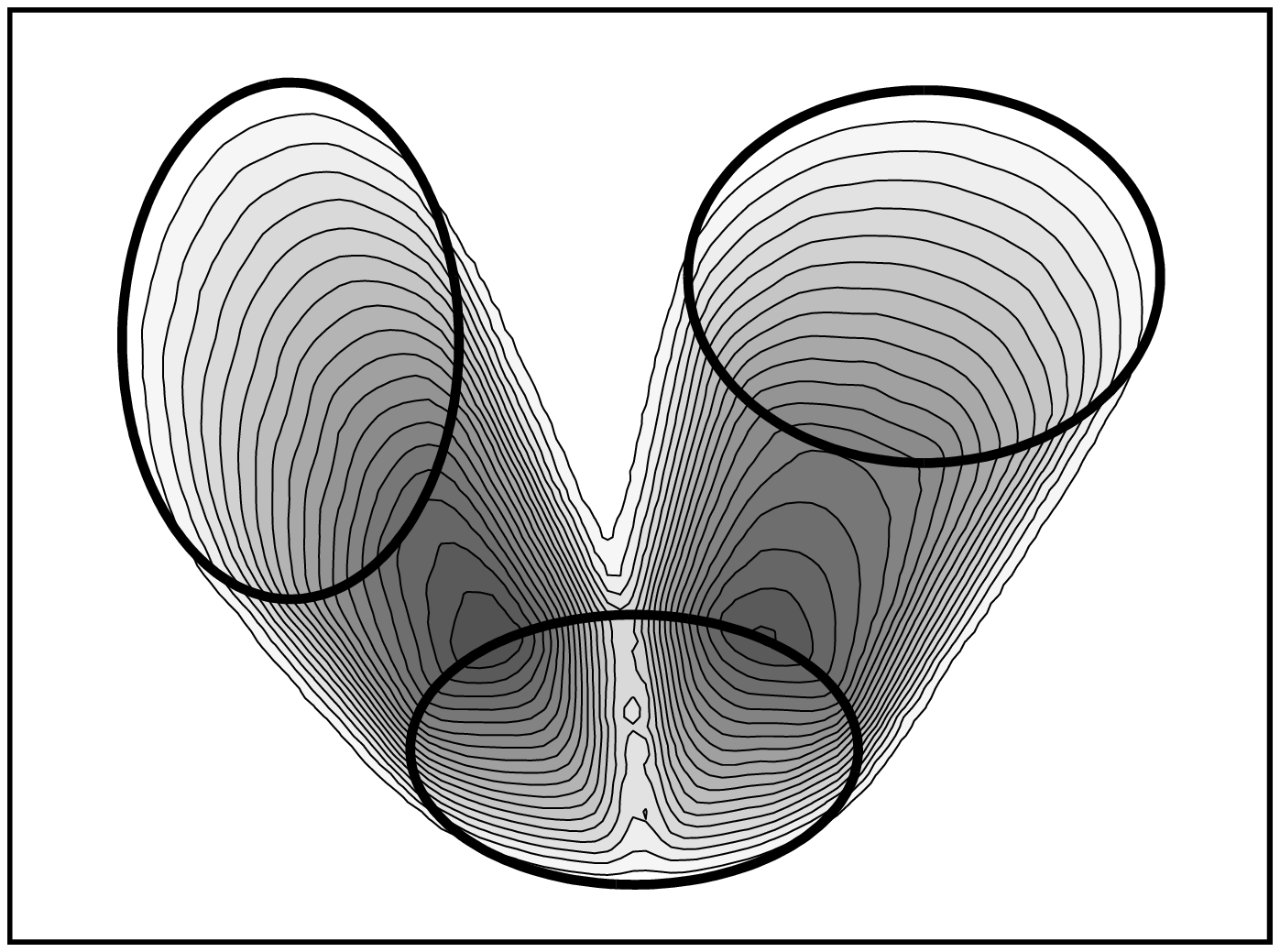}
\end{center}
\caption{Solution to Monge-Kantorovich equations. Left: the
optimal flux and levels of Kantorovich potential. Right: contour
plot of $|\bm{q}|$.}\label{fig1}
\end{figure}

\begin{figure}[htb]
\begin{center}
\includegraphics[width=7.0cm,height=7.0cm]{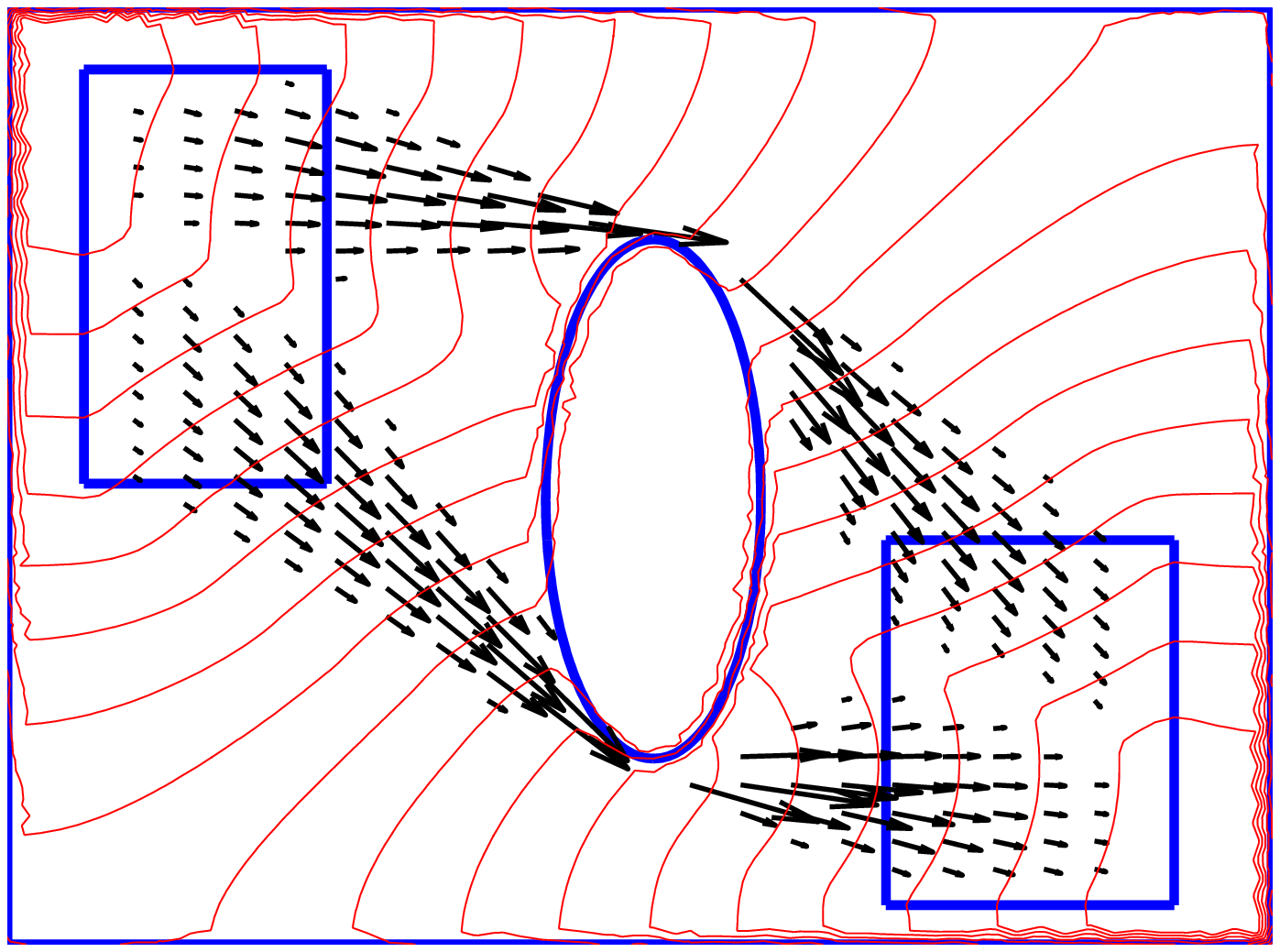}\includegraphics[width=7.0cm,height=7.0cm]{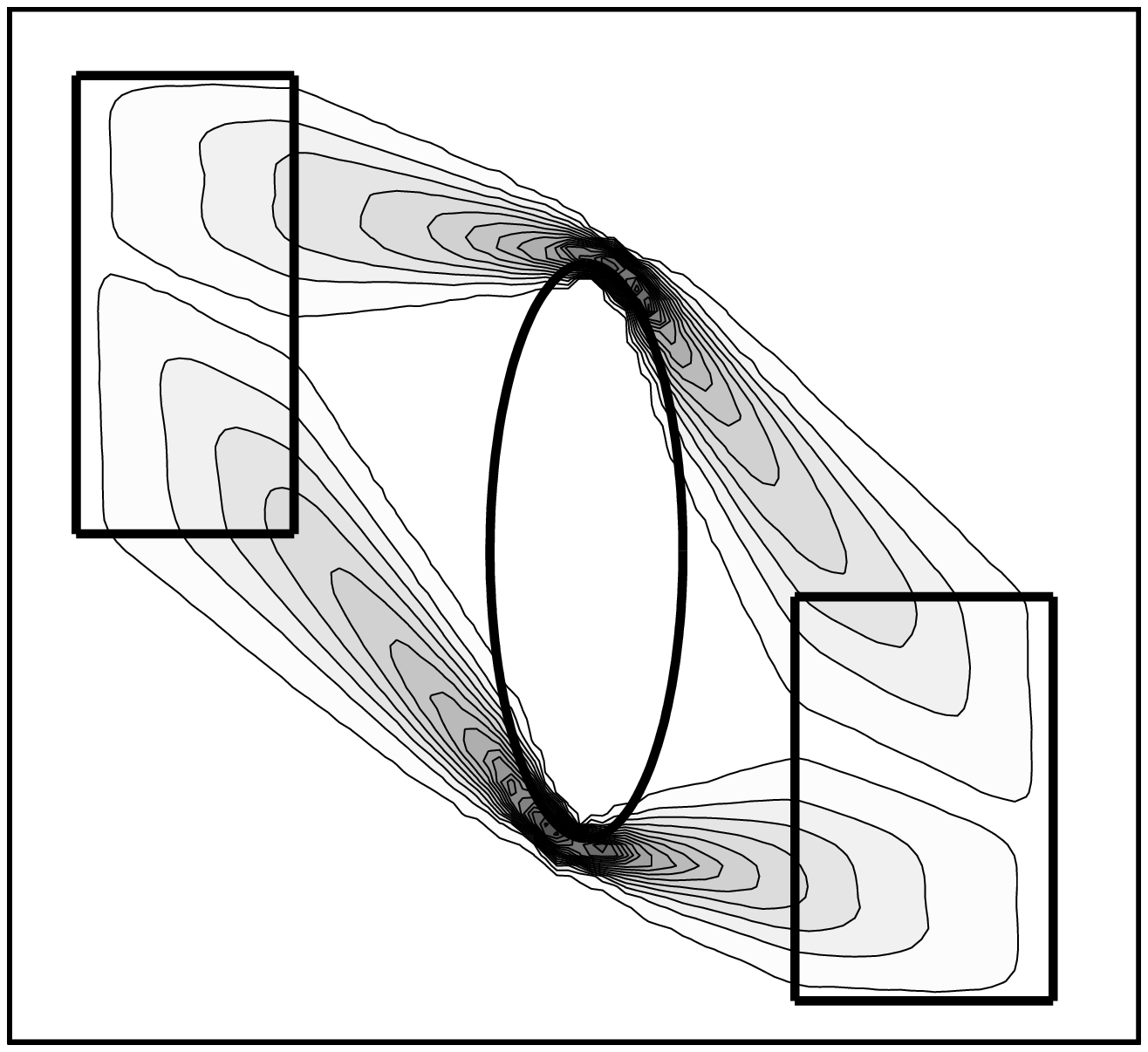}
\end{center}
\caption{As in fig. \ref{fig1}; problem with an obstacle to mass
transportation.}\label{fig2}
\end{figure}

\begin{figure}[htb]
\begin{center}
\includegraphics[width=7.0cm,height=7.0cm]{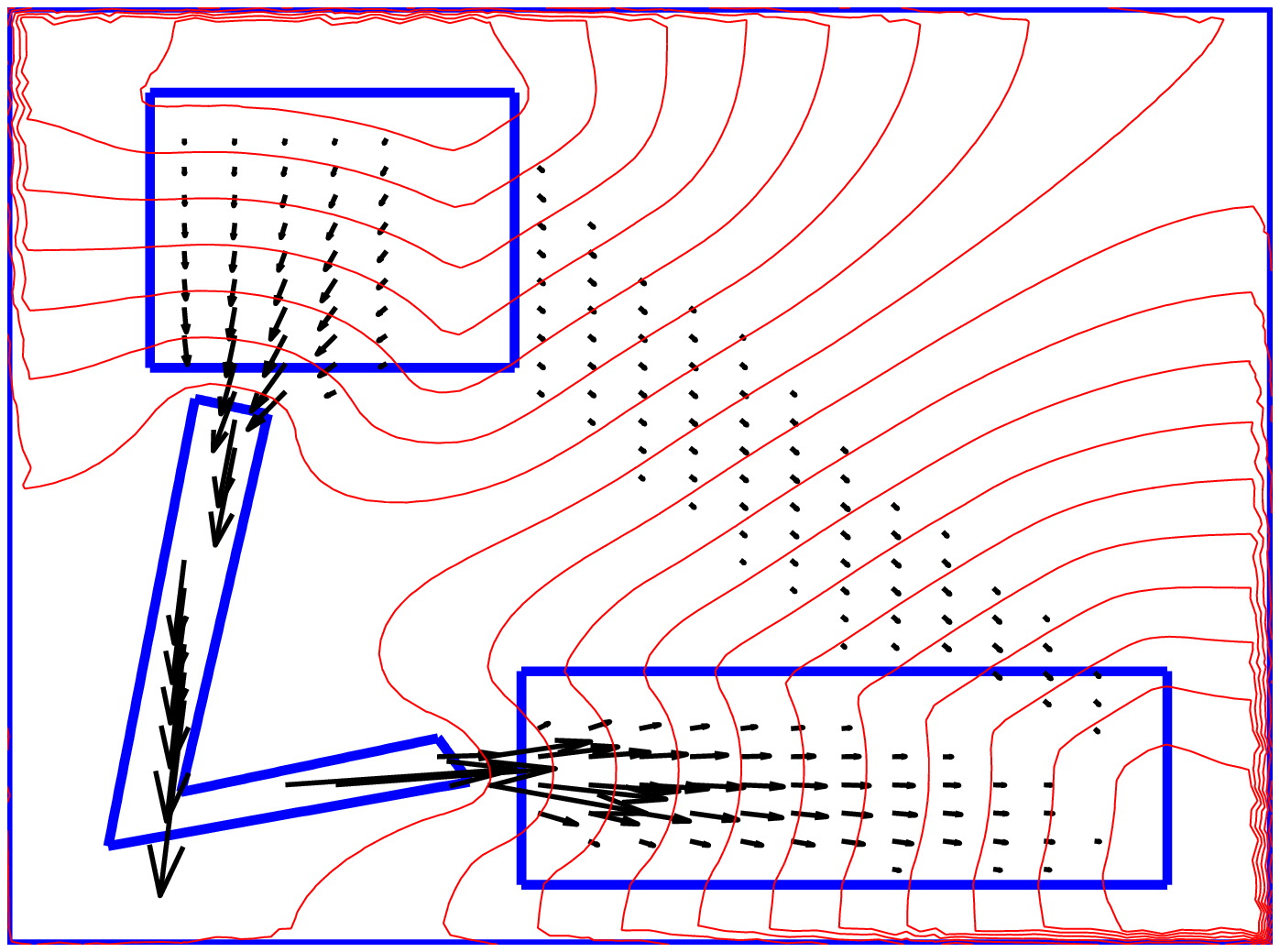}\includegraphics[width=7.0cm,height=7.0cm]{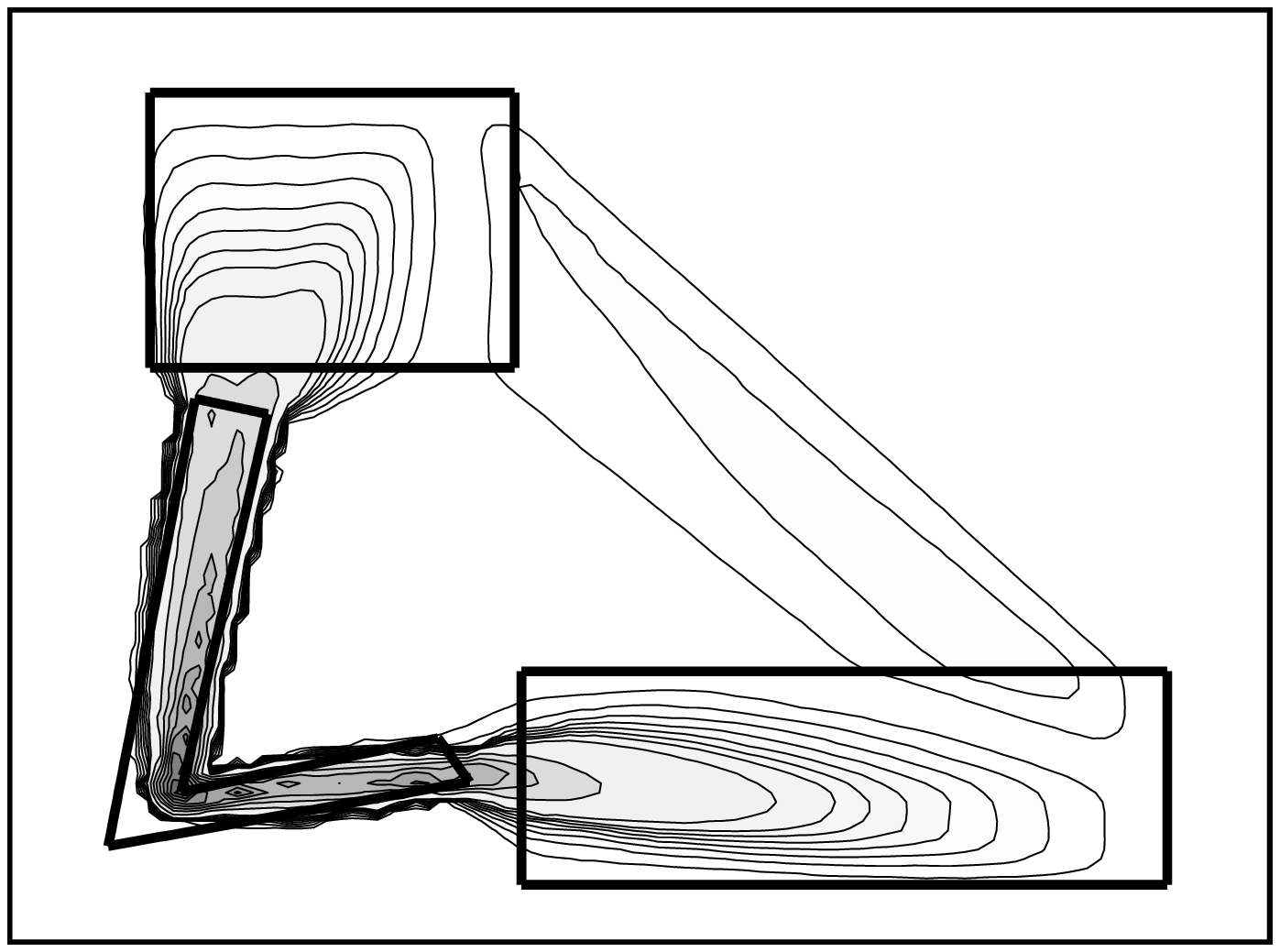}
\end{center}
\caption{As in fig. \ref{fig1}; problem with a region of cheap
transportation.}\label{fig3}
\end{figure}

In all examples below, the computations were performed in a square
surrounded by an obstacle (we set $k=10^6$ in a thin border around
the square). Matlab Partial Differential Equation Toolbox was used
for generation of finite element meshes.
\\
\textit{Example 1}. Support of $f^+$ consists of two upper
ellipses in which the density of $f^+$ is 1 (see Fig. \ref{fig1}).
The mass from $f^+$ is to be uniformly distributed in the third
ellipse below.
\\
\textit{Example 2}. Let $f^+$ and $f^-$ be uniformly distributed
in upper and lower rectangles, respectively, the ellipse between
them is an obstacle (see Fig. \ref{fig2}). In this case the
transport density is not absolutely continuous: it becomes
concentrated on part of the obstacle boundary (see also \cite{F}
for analysis of transport density in a sandpile growth problem
with an obstacle).
\\
\textit{Example 3}. Let $f^+$ and $f^-$ be uniformly distributed
in upper and lower rectangles, respectively, and the third polygon
in Fig. \ref{fig3} be a ``highway" where the transportation cost
is significantly reduced. To model this situation, we set $k=0.01$
in the polygon (which makes transportation through this area
hundred times cheaper). Not surprisingly, the optimal
transportation in this case is concentrated mostly upon the
highway.

\section{Conclusion}
We used formally the similarity of the Monge-Kantorovich equations
and a model for sand surface evolution to present solutions to
Monge-Kantorovich equations as stationary points of a dynamical
system and find them numerically. We note that, although it
remains to prove rigorously the solutions to non-stationary
problem do have the $t\rightarrow\infty$ limit, numerical examples
confirm this statement and demonstrate the efficiency of our
approach. Our hypothesis is the convergence, in fact,  occurs in a
finite time.

It seems interesting to give also an ``optimal mass
transportation" interpretation of sand model equations. Let
$f^+(x,t)$ and $f^-(x,t)$ be, respectively, the production and
consumptions rates of some homogeneous commodity and $u_0(x)$ its
initial distribution. Let, at any time moment, the local price of
this commodity be defined as $-u(x,t)$, i.e., the price is
positive if there is a deficit, $u(x,t)<0$, and is negative if
there is a surplus $u(x,t)>0$ which needs storing.

Whenever the price difference between two arbitrary points, $x$
and $y$, exceeds the cost of transportation of one commodity unit
from one of these points to another, $c(x,y)= |x-y|$, the prices
adjust themselves instantaneously, since a cheaper price is
available by buying elsewhere and transporting. Condition of price
equilibrium can thus be written in the form $|u(x)-u(y)|\leq
|x-y|$ or, locally, $|\nabla u|\leq 1$. Assuming the initial
distribution of goods $u_0$ satisfies this condition, we describe
its further evolution driven by production and consumption as
$$\partial_t u+\nabla\cdot\bm{q}=f^+-f^-,\ \ \ u|_{t=0}=u_0(x)$$
and assume that the mass flux $\bm{q}$ is always directed towards
the price gradient: $\bm{q}=-a\nabla u$, where $a(x,t)\geq 0$ is
an unknown scalar function. Finally, if $|u(x)-u(y)|<|x-y|$, the
transportation from one of these points to another can only
increase the cost and is unprofitable. Locally, this can be
reformulated as the condition $|\nabla u|<1\longrightarrow
a(x,t)=0 $, which brings us to the critical-state model
(\ref{h_eq})-(\ref{h_rel}) in which all transport occurs at the
border of equilibrium.



\end{document}